\documentclass[english,12pt,twoside]{article}
\usepackage{amssymb,amsbsy,amsmath,amsfonts,amssymb,amscd}
\usepackage{latexsym}
\usepackage{euscript}
\usepackage{exscale}
\usepackage{epsfig}
\usepackage{amsthm}
\usepackage[english]{babel}

\newtheorem{theo}{Theorem}[section]

\newtheorem{lem}{Lemma}[subsection]
\newtheorem{prop}{Proposition}[section]

\newtheorem{cor}{Corollary}[subsection]

\newcommand{\bb}{\mathbb}

\newcommand{\St}{\mathrm St}
\newcommand{\StI}{{\mathrm St}_{I}}

\renewcommand{\ker}{\mathrm{Ker}}

\newcommand{\al}{\alpha}
\newcommand{\be}{\beta}

\newcommand{\De}{\Delta}

\newcommand{\om}{\omega}
\newcommand{\Om}{\Omega}

\title{Generalized Steinberg representations of split reductive linear algebraic groups}
\author{Y. A\"it Amrane}

\begin{document}

\maketitle

\selectlanguage{english}
\begin{abstract}
We generalize results of P. Schneider and U. Stuhler for $GL_{l+1}$ to a reductive
 algebraic group G defined and split over a non-archimedean local field $K$. Following
 their lines, we prove that the generalized Steinberg representations of $G$
 with coefficients in any abelian group are cyclic. When $G$ is semi-simple
of adjoint type, we give an expression of these representations, whenever it is possible and
in particular for those that are of maximal degree, in terms of the parahoric subgroups of $G$.
\end{abstract}

\selectlanguage{english}
\begin{center}
{\bf Introduction}
\end{center}

Let $K$ be a non-archimedean local field, i.e. a
finite extension of ${\mathbb Q}_{p}$ or ${\mathbb F}_{p}((t))$.
Let $G$ be the group of $K$-rational points of a connected reductive
linear algebraic group defined and split over $K$. Let $T$ be a split maximal
torus in $G$ and $P$ be a Borel subgroup of $G$ containing $T$. Let $\Phi=\Phi(G,T)$ be
the root system of $G$ with respect to $T$. Denote $\Delta=\{1,2,\ldots ,l\}$, where $l$
is the semi-simple rank of $G$. Let $\Phi_0=\{\alpha_i;\, i\in \Delta\}$ be the basis
of simple roots defined by the Borel $P$. Let $W$ be the Weyl group of $\Phi$. For any
$I\subseteq \Delta$, let $P_I$ be the parabolic subgroup of $G$ generated by $P$ and the
reflections $s_{\alpha_i}$, $i\in I$.

Let $M$ be an abelian group. For any $I\subseteq \Delta$, denote by $C^{^\infty}(G/P_I,\,M)$ the space
of $M$-valued locally constant functions on $G/P_I$. The action of $G$ on $C^{^\infty}(G/P_I,\,M)$ is
induced by its left action on the quotient $G/P_I$. The Steinberg representation of $G$ associated with
the parabolic $P_I$ is the $M[G]$-module :
$$
\StI(M):=\displaystyle\frac{C^{^\infty}(G/P_I,\,M)}{\sum_{i\in
\Delta\backslash I}C^{^\infty}(G/P_{I\cup\{i\}},\, M)}.
$$
When $I=\emptyset$, this is the usual Steinberg representation of $G$.

Let $B$ be the Iwahori subgroup of $G$ which is associated with the Borel subgroup $P$ by reduction
modulo the maximal ideal of the valuation ring of $K$. For any $I\subseteq \Delta$, let $B_I$ be the
parahoric subgroup of $G$ generated by $B$ and the reflections $s_{\alpha_i}$, $i\in I$. Since we are
working with any coefficients $M$, the theory of Bernstein, Borel and Casselman of representations generated
by them Iwahori fixed vectors can not be applied. Instead, P. Schneider and U. Stuhler, in their paper \cite{Schneider1},
for the group $GL_{l+1}(K)$, by exploring the family $\{gB_IP_I/P_I;\, g\in G\}$ they proved that the generalized
Steinberg representation of $GL_{l+1}(K)$ associated with the parabolic $P_I$ is  generated by the characteristic
function $\chi_{B_IP_I}$. In this paper, after establishing some technical results, we show that their method works as well
for our more general group $G$ (Theorem \ref{cyclicI}) :

\begin{theo}
The $M[G]$-module $C^{^\infty}(G/P_{I},\,M)$ is generated by
the characteristic function $\chi_{B_{I}P_{I}}$.
\end{theo}

Let $C_c^{^\infty}(G/B_I,\, M)$ be the space of $M$-valued compactly supported locally constant functions on $G/B_I$.
The action of $G$ on $C_c^{^\infty}(G/B_I,\, M)$ comes from left translations on $G/B_I$. Consider the $M[G]$-homomorphism
$$
\Psi_I :C_c^{^\infty}(G/B_I,\,M) \rightarrow C^{^\infty}(G/P_I,\, M)
$$
which to $\chi_{B_I}$ associates $\chi_{B_IP_I}$. Notice that by the above theorem it is surjective.

Assume that $G$ is semi-simple of adjoint type.  Let $\omega_i$, $i\in \Delta$, be the fundamental
coweights associated with the basis $\Phi_0$ of the root system $\Phi$. For any $i\in \Delta$, let $t_i\in T$
be the element that acts on the Bruhat-Tits building of $G$ as the translation by $\omega_i$. Again, following
the work for $GL_{l+1}(K)$ of P. Schneider and U. Stuhler in [loc.cit.], we compute the kernel of the above homomorphism
to obtain for our group $G$ (Theorem \ref{kernelI}) :

\begin{theo}
Let $I\subseteq \Delta$ be such that any positive root
$\al=\sum_{i\in I}m_i\al_i$ in the root subsystem $\Phi_I$ satisfy
$m_i\leq 1$ for any $i\in I$. The kernel of the surjective map
$$
H_I : C_{c}^{^\infty}(G/B_I,\,M) \longrightarrow
C^{^\infty}(G/P_I,\,M)
$$
is the $M[G]$-submodule $R_I$ of $C_c^{^\infty}(G/B_I,M)$ generated
by the functions $\chi_{Bt_jB_I}-\chi_{B_I}$, $1\leq j \leq l$.
\end{theo}

Under the assumption on the root subsystem $\Phi_I$ in this theorem, we obtain as corollary (Corollary \ref{stI})
an expression of the Steinberg representation of $G$ associated to $P_I$ in terms of the parahoric subgroups of $G$.
Particulary, when $I=\emptyset$, we get (Corollary \ref{st}) :

\begin{cor} We have an $M[G]$-isomorphism
$$
St(M)\cong \frac{C_c^{^\infty}(G/B,M)}{R+\sum_{i\in \Delta}
C_c^{^\infty}(G/B_{\{i\}},M)}
$$
where $R$ is the $M[G]$-submodule of $C_c^{^\infty}(G/B,M)$
generated by the functions $\chi_{Bt_jB}-\chi_{B}$, $1\leq j \leq
l$.
\end{cor}

In the first section, we introduce some notations and the notion of
generalized Steinberg representations of $G$. In the second section,
we recall the geometric definition of the Bruhat-Tits building
associated to $G$ and the main properties we need in this paper. We
establish some technical results essentially about the parabolic and
parahoric subgroups of $G$. We also generalize to $G$ the theorem
1.7 established for $GL_{l+1}(K)$ in \cite{Yacine}, namely that
there is a canonical one-to-one correspondence between the double
classes in the linear Weyl group $W$ and the Iwasawa cells in $G$.
In the third section, we prove that the generalized Steinberg
representations of $G$ are cyclic, and then give under some
conditions the generalized Steinberg representations in terms of
parahoric subgroups.\bigskip

It is a great pleasure to express my deep gratitude to Peter Schneider. He invited me to M\"unster University where
 I had the opportunity to learn more about the theory of algebraic
groups and of the associated Bruhat-Tits buildings. He also gave me very helpful comments on this subject.

\section{Generalized Steinberg representations}\label{GSR}

Let $K$ be a non-archimedean local field. Let $G$ be the group of
$K$-rational points on a reductive linear algebraic group defined
and split over $K$.  Let $l$ be its semi-simple rank.

For a parabolic subgroup $Q$ of $G$, denote by $C^{^\infty}(G/Q,\,
M)$ the space of $M$-valued locally constant functions on $G/Q$. The
action of $G$ on $C^{^\infty}(G/Q,\, M)$ is induced by its left
action on the quotient $G/Q$. Notice that if $Q'$ is a parabolic
subgroup of $G$ that contains $Q$ we have a natural
$M[G]$-monomorphism $C^{^\infty}(G/Q',\, M) \rightarrow
C^{^\infty}(G/Q,\, M)$. The Steinberg representation of $G$
associated with the parabolic subgroup $Q$ is the $M[G]$-module
$$
\St_Q(M)=\displaystyle\frac{C^{^\infty}(G/Q,\,M)}{\sum_{Q'}
C^{^\infty}(G/Q',\, M)}
$$
where the sum is taken over the parabolic subgroups $Q'$ of $G$
containing $Q$. When $Q$ is Borel this is the usual Steinberg representation of $G$.

Let $T$ be a split maximal torus in $G$. Let $X^{*}(T)$,
respectively $X_{*}(T)$, denote the group of characters,
respectively cocharacters, of $T$. By \cite[Ch.III, \S 8.6]{Borel2},
we have a perfect pairing of abelian groups
$$
\langle \cdot ,\cdot\rangle : X_*(T)\times X^*(T)\longrightarrow
{\mathbb Z}
$$
given by $\langle \lambda,\chi \rangle=m$ if $(\chi \circ
\lambda)(x)=x^m$.

Let $V=X_*(T)\otimes {\mathbb R}$ and
identify its dual vector space $V^*$ with $X^*(T)\otimes {\mathbb
R}$.
Let $\Phi=\Phi(G,T)$ be the root system of $G$ with respect to $T$. To every root
$\alpha\in \Phi$ corresponds a coroot $\Check{\alpha} \in V$ and a convolution that
acts on $V$ by :
$$
s_\alpha(x)=x-\langle x,\alpha \rangle \Check{\alpha}.
$$
Let $W=N_G(T)/T$ be the Weyl group of $(G,T)$. We can see that $N_G(T)$ acts on $X_*(T)$ by conjugation and then $T$ acts trivially.
So we have an action of $W$ on $V$ by automorphisms. We can identify $W$ to the Weyl group of $\Phi$.

Fix a Borel subgroup $P$ of $G$ containing $T$. Let
$\Delta=\{1,2,\ldots ,l\}$ and $\Phi_0=\{\al_i;\, i\in \Delta\}$ the
set of the simple roots defined by the Borel $P$. The Weyl group $W$
is generated by the reflections $s_{\alpha_i}$, $i\in \Delta$. For
any $I\subseteq \Delta$, let $W_I$ be the subgroup of $W$ generated
by $s_{\alpha_i},\, i\in I$. We have the following well-known facts
(see \cite{Borel2} or \cite{springer}) :
\begin{enumerate}
\item[(i)] The set
$$
P_I = PW_IP = \coprod_{w\in W_I}PwP
$$
is a parabolic subgroup of $G$ containing $P$.

\item[(ii)] For any root $\alpha \in \Phi$ there exists an isomorphism $x_\alpha$ of the additive group of $K$ onto a unique
closed subgroup $U_\alpha$ of $G$ such that
    \begin{equation}
    t x_\alpha(a) t^{-1}=x_\alpha(\alpha(t)a) \qquad \textrm{for any } a\in K \textrm{ and } t\in
    T.
    \end{equation}\label{txt}
     The group $G$ is generated by $T$ and the $U_\alpha$'s, $\alpha\in \Phi$. A Chevalley system $(x_\alpha)_{\alpha\in \Phi}$
     in $G$ is a choice of such isomorphisms. Also, see \cite[Ex. 8.1.12]{springer}),
for any $w\in W$ and any $\alpha \in \Phi$ there is $c_{w,\alpha}\in
K^*$ such that :
\begin{equation}\label{wxw}
w x_{\alpha}(a) w^{-1}=x_{w\alpha}(c_{w,\alpha}a), \qquad
\textrm{for any } a\in K.
\end{equation}
\item[(iii)] Let $\Phi^+$ (resp. $\Phi^-$) be the set of positive roots (resp. negative roots) of $\Phi$ with respect to $P$.
Denote by $U^+$ (resp. $U^-$) the unipotent radical of the Borel $P$
(resp. of the opposite Borel of $P$). For any ordering on $\Phi$,
the product map induces bijections
\begin{equation}\label{bij}
\prod_{\alpha\in \Phi^{\pm}}U_\alpha \xrightarrow{\sim} U^{\pm}.
\end{equation}
More generally, see \cite[\S\,0.17.(iii)]{Landvogt}, for $I\subseteq
\Delta$, if we denote by $U_I^+$ (resp. $U_I^-$) the unipotent
radical of the parabolic  $P_I$ (resp. of the opposite parabolic),
the product map induces bijections
\begin{equation}\label{bijI}
\prod_{\alpha\in \Psi_I^{\pm}} U_\alpha \xrightarrow{\sim} U_I^{\pm}
\end{equation}
where $\Psi_I^{\pm}=\Phi^{\pm}-\Phi_I$ and $\Phi_I$ is the root subsystem of $\Phi$ of the roots that are linear combinations
of the $\alpha_i$'s, $i\in I$.
\item[(iv)] Any parabolic subgroup of $G$ containing $P$ is some $P_I$ for a unique $I\subseteq \Delta$. Any parabolic
subgroup of $G$ is conjugate to some $P_I$.
\end{enumerate}

Since of this last fact, from now on, we need only to consider the
Steinberg representations of $G$ associated with the $P_I$'s,
$I\subseteq \De$, they will be denoted by :
$$
\StI(M):=\St_{P_I}(M)=\displaystyle\frac{C^{^\infty}(G/P_I,\,M)}{\sum_{i\in
\Delta\backslash I}C^{^\infty}(G/P_{I\cup\{i\}},\, M)}.
$$

\section{Bruhat-Tits buildings}

Let $K$ be a complete field with respect to a discrete valuation
$\om$ such that $\om(K^*)=\mathbb Z$. Let $O$ be its valuation ring,
$\pi\in O$ an uniformizing element and $k=O/\pi$ the residue field.

Let $G$ be the group of $K$-rational points of a reductive
linear algebraic group defined and split over $K$. Only for simplicity,
assume $G$ is semi-simple of rank $l$.
Let $T$ be a split maximal torus in $G$ and $N=N_G(T)$ its normalizer in $G$.

For complete presentations on the Bruhat-Tits building attached to a reductive group, the reader may refer to the original
construction of F. Bruhat and J. Tits, \cite{Bruhat1} and \cite{Bruhat2}. For more concise presentations the reader may refer
to \cite{Landvogt} and \cite{Schneider2}.

\subsection{The building}\label{building}

Let $A$ be an affine space under the vector space $V$. It is a Coxeter complex, where the walls of the chambers are the hyperplanes :
$$
H_{\al,r}:=\{x\in A;\; \langle x,\al \rangle +r=0\},\quad \al\in
\Phi, r\in {\mathbb Z}.
$$
There is a unique group homomorphism $\nu : T \rightarrow
X_*(T)\subset V$ such that :
\begin{equation}\label{TA}
\langle \nu(t),\chi \rangle=-\om(\chi(t))
\end{equation}
for any $t\in T$ and any $\chi\in X^*(T)$. The torus $T$ acts on $A$ by translations :
\begin{equation}\label{xnut}
tx:=x+\nu(t); \quad t\in T, x\in A.
\end{equation}
This action can be extended to an action of $N$ on $A$ by affine
automorphisms. The space $A$ equipped with its simplicial structure
and of the affine action of $N$ is called the fundamental apartment
of the building.

For any nonempty subset $\Om\subseteq A$, let $f_{\Om}: \Phi
\rightarrow {\bb R}\cup\{\infty\}$ be the map defined by
$$
f_{\Om}(\al)=-\inf_{x\in \Om}\{\al(x)\}.
$$
For any $\alpha\in  \Phi$ and any $r\in \mathbb R$, denote $
U_{\alpha,r}=\{x_{\alpha}(a);\, a\in K\textrm{ and }\omega(a)\geq
r\} $ and denote  by $U_{\Om}$ the subgroup of $G$ generated by the
subgroups $U_{\al,f_{\Om}(\al)}$ for $\al \in \Phi$.

The Bruhat-Tits building $X$ of $G$ is obtained by taking the
quotient of the direct product $G\times A$ by the equivalence
relation :
$$
(g,x)\sim (h,y) \quad \textrm{if there is} \; n\in N
\;\;\textrm{s.t.} \quad y=nx\;\; \textrm{and}\;\; g^{-1}hn\in U_{\{x\}}.
$$
The action of $G$ on X is given by
$$
g \cdot \textrm{Class}(h,y)=\textrm{Class}(gh,y) \quad \textrm{for }g\in G \textrm{ and } (h,y)\in G\times A.
$$
The apartment $A$ is viewed as a subset of $X$ by the embedding
$x\mapsto \textrm{Class}(1,x)$. The other apartments of $X$ are the
$gA$'s, $g\in G$.

Let $\Omega \subseteq X$ be a non-empty subset. We define ${\mathcal P}_\Omega$ to be the pointwise stabilizer of $\Omega$ in $G$
$$
{\mathcal P}_\Omega :=\{g\in G;\, gx=x \textrm{ for any }x\in \Omega\}.
$$
In case $\Omega\subseteq A$, by the action of $N$ on the apartment
$A$, we can also define $N_\Omega$ to be the pointwise stabilizer of
$\Omega$ in $N$ and we have ${\mathcal P}_\Omega =U_\Omega
N_\Omega$.

Recall from \S\,\ref{GSR} that $\Delta=\{1,2,\ldots l\}$ and that
$\Phi_0=\{\alpha_i;\, i\in \Delta\}$ is the set of the simple roots
corresponding to the Borel $P$, that $\Phi=\Phi^+\cup \Phi^-$ is a
decomposition of $\Phi$ into positive and negative roots with
respect to the fixed Borel $P$, and that $U^{+}$ (resp. $U^-$) is
the unipotent radical of $P$ (resp. of the opposite Borel of $P$).
Denote $U_\Omega^{\pm}=U_\Omega \cap U^{\pm}$. For any ordering of
the factors the product map induces bijections
\begin{equation}\label{rootdecompositionUomega}
\prod_{\alpha\in \Phi^{\pm}}U_{\alpha,f_\Omega(\alpha)} \rightarrow
U_\Omega^{\pm},
\end{equation}
and the group ${\mathcal P}_\Omega$ has the following Iwahori decomposition :
\begin{equation}\label{IwahoriOmega}
{\mathcal P}_\Omega =U_\Omega^{-}U_\Omega^{+}N_\Omega.
\end{equation}

Since the action of $G$ on the building is in fact transitive on the chambers (simplices of maximal dimension),
we will focus only on what we call the fundamental chamber. Let $\tilde{\alpha}$ be the longest root in $\Phi$.
The fundamental chamber $C$ is the chamber of the apartment $A$ defined by
$$
C:=\{x\in A;\; \langle x, \alpha_i\rangle \gneq 0 \textrm{ for any
}i\in \Delta \textrm{ and } \langle x,\tilde{\alpha}\rangle \lneq
1\}.
$$
Also, for any $I\subseteq \Delta$, define $C_I$ to be simplex which is the face of $C$ given by
$$
C_{I}=\{x\in C;\, \langle x,\al_{i}\rangle =0 \textrm{ for }i\in I\}.
$$
Notice that for $\Omega=C$, it is not hard to see that when
$\alpha\in \phi^+$ we have $f_\Omega(\alpha)=0$, and when $\alpha\in
\Phi^-$ we have $0\lneq f_\Omega \leq 1$ and then
$U_{\alpha,f_\Omega(\alpha)}=U_{\alpha,1}$. Therefore, for any
ordering on $\Phi^+$ and on $\Phi-$, by
(\ref{rootdecompositionUomega}) the product map induces bijections :
\begin{equation}\label{rootdecompositionUC}
U_C^-=\prod_{\alpha\in \Phi^-}U_{\alpha,1} \quad \textrm{ and }\quad
U_C^+=\prod_{\alpha\in\Phi^+}U_{\alpha,0}.
\end{equation}

According to Casselman \cite[Proof of Prop. 1.4.4.]{Casselman},
there exists a smooth group scheme ${\bf G}$ over $O$ such that
$G\cong {\bf G} \times \textrm{Spec}(K)$. For each $n\geq 0$, let
\begin{equation}\label{reductionmodpin}
{\bf G}(O) \rightarrow {\bf G}(O/\pi^{n}O)
\end{equation}
be the reduction homomorphism modulo $\pi^{n}$. For any $I\subseteq
\Delta$, let $B_{I}$ be the inverse image by the reduction mod $\pi$
of the parabolic $P_I(k)$. In fact $B_{I}$ is nothing else than the
pointwise stabilizer ${\mathcal P}_{C_I}$ in $G$ of the simplex
$C_I$.

For later use, we need also to consider what we call  the fundamental vectorial chamber
$$
{\mathfrak C}=\{x\in A;\,\langle x,\al_{i}\rangle \gneq 0 \textrm{ for any } i\in
\Delta\},
$$
and its faces which are the conical cells ${\mathfrak C}_{I}$,
$I\subseteq \Delta$, given by
$$
{\mathfrak C}_{I}=\{x\in A;\,\langle x,\al_{i}\rangle =0 \textrm{
for }i\in I\;\textrm{and } \gneq 0\textrm{ for } i\in \De-I\}.
$$
Notice that $P_I$ is the stabilizer of the simplex $({\mathfrak C}_{I})_{\infty}$
defined at infinity by the conical cell ${\mathfrak C}_I$, see \cite{Borel3}. Therefore, see
 \cite{Yacine0}, the intersection $B_{I}\cap P_{I}$ is the pointwise stabilizer in $G$  of the conical cell
${\mathfrak C}_{I}$, i.e.
$$
{\mathcal P}_{{\mathfrak C}_{I}}=B_{I}\cap P_{I}.
$$

To finish this paragraph, recall also that we have the following important decompositions in $G$ : \medskip\\
{\bf Bruhat-Tits decomposition :} For any $I\subseteq \Delta$, we
have
$$
B_I=\coprod_{w\in W_I}BwB. \medskip
$$
{\bf Cartan decomposition :}
Consider the monoid
$$
T^{++}=\{t\in T;\, \om \circ \al_{i}(t) \leq 0\; \textrm{for any}\;
i\in \Delta\}.
$$
We have
$$
G=\coprod_{t\in T^{++}}={\bf G}(O)t {\bf G}(O).
$$

\subsection{Iwahori decomposition}

Recall form \S\,\ref{GSR} that for any $I\subseteq \De$, we have denoted :\smallskip\\
$\bullet$ $\Phi_{I}$ the set of all roots in $\Phi$ that are linear combinations
of the simple roots $\al_{i}$, $i\in I$.\smallskip\\
$\bullet$ $\Psi_{I}^+=\Phi^+-\Phi_{I}$ and $\Psi_{I}^-=\Phi^--\Phi_{I}$.\smallskip\\
$\bullet$ $U^{+}_{I}$ (resp. $U^-_{I}$) the unipotent group generated by the $U_{\al}$'s,
$\al\in \Psi_{I}^+$  (resp. $\al \in \Psi_{I}^-$).

Let $I\subseteq \Delta$. It is essentially proven, see \cite[Proof
of Th. 2.5]{Iwahori}, that the Iwahori subgroup $B$ has an Iwahori
factorization with respect to the parabolic $P_{I}$, i.e. we have
\begin{equation}\label{IwahoriI}
B=(B\cap U_{I}^-)(B\cap P_{I}).
\end{equation}
Since $B={\mathcal P}_C$, we have
\begin{equation}\label{ChevalleyI}
B\cap U_I^{-} =U_\Omega \cap U_I^-=\prod_{\al \in
\Psi_I^{-}}U_{\al,1},
\end{equation}
and using formulas (\ref{txt}) and (\ref{TA}), we can see that we
have
\begin{equation}\label{tUI}
t(B\cap U_{I}^-)t^{-1} \subseteq B\cap U_{I}^{-} \quad \textrm{for any}\; t\in T^{++}.
\end{equation}

\begin{lem}\label{wbuwb}
For any $w\in W$ we have $w(B\cap U^-)w^{-1}\subset U_C \subset B$.
\end{lem}
\proof Let $\al\in\Phi^-$ and let $w\in W$. We have
$wU_{\al,1}w^{-1}=U_{w\al,1}$. If $w\alpha\in \Phi^-$, then
$U_{w\al,1}\subseteq U_C^-\subseteq U_C$. If $w\al\in\Phi^+$, then
$U_{w\al,1}\subseteq U_{w\al,0}\subseteq U_C^+\subseteq U_C$. So in
either case, $wU_{\al,1}w^{-1}\subseteq U_C$. By (\ref{ChevalleyI}),
we get
$$
w(B\cap U^-)w^{-1}= \prod_{\al\in \Phi^-}wU_{\al,1}w^{-1}\subseteq U_C\subseteq B.
$$
\qed

\begin{lem}\label{wbuiiwbui}
Let $I\subseteq \De$. For any  $w\in W_{I}$, we have $w(B\cap U_{I}^-)w^{-1}
=B\cap U_{I}^{-}$.
\end{lem}
\proof Let $\al \in \Psi_{I}^-$, then $\al = \sum_{i\in
\De}n_{i}\al_{i}$ with $n_{j}<0$ for some $j\in \De-I$. Therefore $s_{i}\al \in \Psi_{I}^-$ for
any $i\in I$. Then $W_{I}$ stabilize $\Psi_{I}^-$.  \smallskip\\
Let $w\in W_{I}$, let $\al\in \Psi_{I}^-$ and $u\in U_{\al}$. Let
$a\in K$ such that $x_{\al}(a)=u$. By (\ref{wxw}) there is a
constant $c_{w,\al}\in K$ such that
$wx_{\al}(a)w^{-1}=x_{w\al}(c_{w,\al}a)$; thus
$wuw^{-1}=wx_{\al}(a)w^{-1}$ is in $U_{w\al}$ with $w\al \in
\Psi_{I}^-$. Henceforth, since $U_{I}^{-}$ is generated by the
$U_{\al}$'s, $\al\in \Psi_{I}^{-}$, we deduce that $W_{I}$
normalizes $U_{I}^-$. Finally, since we have $U_{I}^- \subseteq U^-$
we get
$$
w(B\cap U_{I}^- )w^{-1}=w(B\cap U^- \cap U_{I}^- )w^{-1}=w(B\cap U^- )w^{-1}
\cap wU_{I}^- w^{-1}
$$
for any $w\in W_{I}$. By Lemma \ref{wbuwb} above we have $w(B\cap
U^-)w^{-1}\subset B$ for any $w\in W$ and so for any $w\in W_{I}$,
and we just established that $W_{I}$ normalizes $U_{I}^-$. \qed

The following proposition is proved for the group $GL_{l+1}$ in
\cite[Lemma 14.i.]{Schneider1}. With the same arguments and using
Lemma \ref{wbuwb} above, we get
\begin{prop}\label{bibuibipi}
 Let $I\subseteq \Delta$. We have the factorization
$$
B_{I}=(B\cap U_{I}^-)(B_{I}\cap P_{I}).
$$
\end{prop}
\proof Notice that for any $w\in W$, by the Iwahori factorization of
$B$ with respect to $P$ and using Lemma \ref{wbuwb} we get
$$
BwB=Bw(B\cap U^-)(B\cap P)=Bw(B\cap P).
$$
By the Bruhat-Tits decomposition and the Iwahori factorization of
$B$ with respect to $P_I$, we then get
$$
B_I=\coprod_{w\in W_I}BwB=\coprod_{w\in W_I}(B\cap U_I^-)(B\cap
P_I)w(B\cap P).
$$
Since for any $w\in W_I$ we have $(B\cap P_I)w(B\cap P)\subseteq
B_I\cap P_I$, we deduce that $B_{I}=(B\cap U_{I}^-)(B_{I}\cap
P_{I})$.

\qed

Moreover, if for any $I\subseteq \Delta$ we define
$$
T_{I}^{++} =\{t\in T^{++}\,;\; \om \circ \al_{i}(t)=0 \;\textrm{for}\;i\in I\},
$$
we have
\begin{equation}\label{TIBIPI}
t^{-1}(B_I \cap P_I)t \subseteq B_I\cap P_I \qquad \textrm{for any } t\in T_{I}^{++}.
\end{equation}
Indeed, this comes from the fact that $B_{I}\cap P_{I}$ is the
stabilizer of the conical cell ${\mathfrak C}_{I}$ and geometrically
we can check easily that for any $t\in T_I^{++}$ the translation of
${\mathfrak C}_I$ by $\nu(t)$ is contained in ${\mathfrak C}_I$.

\subsection{Generalized Iwasawa decomposition}

In this paragraph we shall follow the principal steps in \cite[\S\,1.2.2.]{Yacine} and use the Iwasawa decomposition of $G$
$$
G=\coprod_{w\in W}BwP
$$
to prove Theorem \ref{iwasawaI} which gives such a decomposition of
$G$ but with respect to $B_I$ on the left and to $P_I$ on the right
for any $I\subseteq \Delta$. All is based on the following result.

\begin{prop}\label{wps}
Let $w\in W$, $\alpha \in \Phi_0$. We have
\begin{enumerate}
\item $wPs_\alpha \subseteq BwP \cup Bws_\alpha P$
\item $s_\alpha Bw \subseteq BwP \cup Bs_\alpha wP$
\end{enumerate}
\end{prop}

\proof

Let $w\in W$, $\al \in \Phi_{0}$. Proving 1. and 2. comes to prove respectively
\begin{enumerate}
\item[1'.] $Ps_{\al}\subseteq B'P \cup B's_{\al}P$ where we set
$B'=w^{-1}Bw$.
\item[2'.] $s_{\al}B \subseteq BP' \cup Bs_{\al}P'$ where we set
$P'=wPw^{-1}$.
\end{enumerate}
For $\al\in \Phi$ we denote by $G_{\al}$ the centralizer in $G$ of
the connected component of $\ker\, \al$. It is a reductive group of
semi-simple rank one, and is generated by $T$, $U_{\al}$ and
$U_{-\al}$. We know, see for example \cite[Th. 6.4.7
(ii)]{springer}, that the intersection of a Borel subgroup $G$
containing $T$ with $G_\alpha$ is a Borel subgroup of $G_\alpha$.

\begin{itemize}
\item {\it Proof of 1'.} Let us prove that we have $Ps_{\al}\subseteq G_{\al}P$.
Give an ordering of $\Phi^+$ such that we have $U^+=\prod_{\be\in
\Phi^+}U_{\be}=U_{\al}\prod_{\be\in
\Phi^+\backslash\{\al\}}U_{\be}$, we then have
$$
U^+s_{\al}=s_{\al}(s_{\al}U_{\al}s_{\al})\prod_{\be\in \Phi^+
\backslash\{\al\}}(s_{\al}U_{\be}s_{\al})=s_{\al}U_{s_{\al}(\al)}
\prod_{\be\in \Phi^+\backslash\{\al\}}U_{s_{\al}(\be)}.
$$
Since we have $s_\alpha(\alpha)=-\alpha$, and by \cite[Ch. VI,\S\,1.6, Cor. 1]{Bourbaki} the reflection  $s_\alpha$
permutes the positive roots in $\Phi^+\backslash\{\alpha\}$, we get
$$
U^+s_{\alpha}=s_{\alpha}U_{-\alpha}\prod_{\be\in
s_{\al}(\Phi^+\backslash\{\al\})}U_{\be}\subseteq
s_{\alpha}U_{-\alpha}U^+ \subseteq G_{\alpha}P.
$$
Thus $Ps_{\al}\subseteq G_{\al}P$. Now, to prove 1', it is enough to
prove
$$
G_{\al}\subseteq (B'\cap G_{\al})(P\cap G_{\al}) \cup (B'\cap
G_{\al})s_{\al}(P\cap G_{\al}).
$$
But, the group $B'\cap G_\alpha$ is the Iwahori subgroup of $G_\alpha$ corresponding to the Borel
$w^{-1}Pw\cap G_{\alpha}$. So, depending on $w$ and $\alpha$, it could be either $B\cap G_\alpha$ which
corresponds to the Borel $P\cap G_{\alpha}$, either $s_\alpha (B\cap G_{\alpha})s_\alpha$ which corresponds
to the opposite Borel of $P\cap G_\alpha$. In either cases, the inclusion above which is an equality is nothing
else than the Iwasawa decomposition in $G_\alpha$.

\item {\it Proof of 2'.} Similarly, let us first prove that we have
the inclusion $s_{\al}B\subseteq BG_{\al}$. By the bijections
(\ref{rootdecompositionUC}), fix an ordering of $\Phi^-$ and of
$\Phi^+$  such that we have
$$
U_C^-=(\prod_{\be\in \Phi^-\backslash\{-\al\}}U_{\be,1})U_{-\al,1}
\qquad \textrm{and}\ qquad  U_C^+=U_{\alpha,0}(\prod_{\be\in
\Phi^+\backslash\{\al\}}U_{\be,0}).
$$
Therefore, as we did in 1' above, we have
$$
s_{\al}U_C^-U_C^+=(\prod_{\be\in
\Phi^-\backslash\{-\al\}}U_{s_\alpha(\be),1})U_{\al,1}(\prod_{\be\in
\Phi^+\backslash\{\al\}}U_{s_\alpha(\be),0})U_{-\al,0}s_{\al}.
$$
Since $s_\alpha$ permutes the negative roots in
$\Phi^-\backslash\{-\alpha\}$ and respectively the positive roots in
$\Phi^+\backslash\{\alpha\}$, we get
$$
(\prod_{\be\in
\Phi^-\backslash\{-\al\}}U_{s_\alpha(\be),1})\subseteq U_C^- \subseteq B \quad \textrm{ resp. } \quad
U_{\al,1}(\prod_{\be\in
\Phi^+\backslash\{\al\}}U_{s_\alpha(\be),0})\subseteq U_C^+ \subseteq B.
$$
Since we also have $U_{-\al,0}s_{\al}\subseteq G_\alpha$, we deduce
$s_{\alpha}U_C^-U_C^+\subseteq BG_{\alpha}$. Therefore
$s_{\al}B\subseteq BG_{\al}$. Now, for the proof of 2', it is enough
to see that we have
$$
G_{\al}\subseteq (B\cap G_{\al})(P'\cap G_{\al}) \cup (B\cap
G_{\al})s_{\al}(P'\cap G_{\al}).
$$
Indeed, as for 1' above, this is again the Iwasawa decomposition in
$G_{\al}$.\qed
\end{itemize}

\begin{cor} Let $w\in W$, $u_1 , \ldots ,u_d \in S$. We have :
\begin{enumerate}
\item $wPu_{1}\cdots u_{d} \subseteq \bigcup_{(l_{1},\ldots ,l_{p})}
Bwu_{l_{1}}\cdots u_{l_{d}}P$
\item $u_{1}\cdots u_{d}Bw \subseteq \bigcup_{(l_{1},\ldots ,l_{p})}
Bu_{l_{1}}\cdots u_{l_{d}}wP$
\end{enumerate}
where  $(l_{1},\ldots ,l_{p})$ runs through all the increasing sequences
of the interval $\{1,\ldots ,d\}$.
\end{cor}
\proof Induct on $d$ and use Proposition \ref{wps}. \qed

\begin{cor}\label{bipibpi}
Let $I_{1},I_{2}$ be subsets of $\De$. For any $w\in W$ we have
$B_{I_{1}}wP_{I_{2}}=BW_{I_{1}}wW_{I_{2}}P$. In particular, for any
$I \subseteq\De$, we have $B_{I}P_{I}=BP_{I}=B_{I}P$.
\end{cor}
\proof See \cite[Cor. 1.6 and Rem. 1.8]{Yacine}. \qed

\begin{theo}\label{iwasawaI}
Let $I_{1},I_{2}$ be subsets of $\Delta$. The map $W\rightarrow
B_{I_1}\backslash G/P_{I_2}$ which to $w$ associates the double
coset $B_{I_1}wP_{I_2}$ induces a bijection :
$$
 W_{I_1}\backslash W/W_{I_2}\sim  B_{I_1}\backslash G/P_{I_2}.
$$

\end{theo}
\proof The proof is the same as in \cite[Th. 1.7]{Yacine}. \qed

\section{Generalized Steinberg representations and parahoric subgroups}

\subsection{Generalized Steinberg representations}

Let $G$ be the group of $K$-rational points of a reductive linear
algebraic group defined and split over $K$ and of semi-simple rank
$l$. In this paragraph, since the parabolic subgroups contain the
unipotent radical of $G$, we can assume $G$ to be semi-simple. By
using central isogenies, we can even assume that $G$ is of adjoint
type.

The following result is established for the group $GL_{l+1}$ in
\cite{Schneider1}, Propositions 7 and 8'.i. Since it was stated for
any $I\subseteq \Delta$ but the proof given in details only in the
case $I=\emptyset$, for the convenience of the reader we will give
the proof here. We can see that it remains valid for our more
general group $G$.

\begin{prop}\label{disjoint}
Let $b,b'\in G(O)$, let $t\in T_{I}^{++}$. If the sets $btB_IP_I$ and $b'tB_IP_I$ are not disjoint then they are equal and we have $btB_I=b'tB_I$ and $bB_I=b'B_I$.
\end{prop}
\proof It suffices to assume the sets $tB_IP_I$ and $btB_IP_I$ to be
not disjoint. By Corollary \ref{bipibpi} and the Iwahori
factorization of $B$ with respect to the parabolic $P_I$, we have
$$
B_IP_I=btBP_I=(B\cap U_I^-)(B\cap P_I)P_I=(B\cap U_I^-)P_I,
$$
and since $t\in T_I^{++}\subseteq P_I$, therefore $P_I=t^{-1}P_I$,
we get
\begin{equation}\label{f1}
tB_IP_I=t(B\cap U_I^-)t^{-1}P_I
\end{equation}
Since we have $t(B\cap U_I^-)t^{-1}\subseteq B\cap U_I^-$, we deduce
that $tB_IP_I\subseteq B_IP_I$ and $btB_IP_I\ \subseteq bB_IP_I$.
Therefore our assumption implies that the sets $B_IP_I$ and
$bB_IP_I$ are not disjoint. Let then $b_0,b_1\in B_I$ and $p\in P_I$
such that $bb_0=b_1p$. Then $p=b_1^{-1}bb_0\in G(O)\cap P_I=B_I\cap
P_I$. Hence $b=b_1pb_0^{-1}\in B_I$. From Proposition
\ref{bibuibipi}, there is $b'\in B\cap U_I^-$ such that $b^{-1}b'\in
B_I\cap P_I$. Now, since $t\in T_I^{++}$, by (\ref{TIBIPI}) we have
$t^{-1}(B_I\cap P_{I})t\subseteq B_I\cap P_I$, and then
$t^{-1}b^{-1}b't\in B_I\cap P_I$. Therefore
\begin{equation}\label{f2}
btB_IP_I=b'tB_IP_I=b't(B\cap U_I^-)t^{-1}P_I.
\end{equation}
From (\ref{f1}) and (\ref{f2}), our assumption implies
\begin{equation}\label{f3}
t(B\cap U_I^-)t^{-1}P_I \cap b't(B\cap U_I^-)t^{-1}P_I \neq
\emptyset.
\end{equation}
Since $b'\in B\cap U_I^-$ and we have $t(B\cap U_I^-)t^{-1}\subseteq B\cap U_I^-$, we must have
$$
t(B\cap U_I^-)t^{-1} \cap b't(B\cap U_I^-)t^{-1} \neq \emptyset.
$$
Thus $b'\in t(B\cap U_I^-)t^{-1}$, and hence $b\in t(B\cap U_I^-)t^{-1}(B_I\cap P_I)$. The assertions
should follow from this.   \qed

Since we assumed $G$ of adjoint type, we know that the natural
inclusion $X_{*}(T)\subseteq P(\Check{\Phi})$ of the group of
cocharacters of T in the coweight lattice  of the root system $\Phi$
is in fact an equality. Since $T$ is split the homomorphism $\nu
:T\rightarrow X_{*}(T)$ given in \S\,\ref{building} by the formula
(\ref{TA}) is surjective. If we denote by $\omega_{j}$ the
fundamental coweights, for any $j\in \Delta$ we can take $t_{j}\in
T$ such that $\nu(t_{j})=\om_{j}$. Also, for any $I\subseteq
\Delta$, put
$$
t_{I}:=\prod_{j\in \De\backslash I}t_{j}.
$$
The following result is also stated in \cite{Schneider1} for the
group $GL_{l+1}(K)$. We will follow the lines of the proof given
there, but still we have to prove some arguments in the more general
situation of our group.
\begin{prop}
The sets $t_I^n B_IP_I/P_I$, $n\geq 0$, form a fundamental system of
neighborhoods of the trivial coset in $G/P_I$.
\end{prop}
\proof For any integer $n\geq 0$, let $K_n$ be the kernel of the
homomorphism reduction modulo $\pi^n$, see \S,\ref{building}.
According to \cite[Prop. 1.4.4]{Casselman} the family of the sets
$K^{(n)}$, $n\geq 0$, is a fundamental system of neighborhoods of
the identity in $G$. Therefore, so is the family of the sets
$B^{(n)}:=B\cap K^{(n)}$, $n\geq 0$. Now, we need only to show that
for any $n$ we have $t_{I}^{n}B_{I}P_{I} \subseteq B^{(n)}P_{I}$. By
the Iwahori decomposition $B_I=(B\cap U_I^-)(B\cap P_I)$ and since
$t_I\in P_I$ we have
$$
t_{I}^{n}B_{I}P_{I}=t_{I}^{n}(B\cap U_{I}^-)t_{I}^{-n}P_{I}.
$$
Let us now prove that $t_I^n(B\cap U_I^-)t_I^{-n}\subseteq B^{(n)}$.
We have $B\cap U_I^-=\prod_{\alpha\in \Psi_I^-}U_{\alpha,1}$, thus
$$
t_I^n(B\cap U_I^-)t_I^{-n}=\prod_{\alpha\in \Psi_I^-}t_I^n
U_{\alpha,1} t_I^{-n}.
$$
Let $u\in U_{\alpha,1}$, thus $u=x_\alpha(a)$ with $a\in K$ and
$\omega(a)\geq 1$. We want to show that $t_I^n u t_I^{-n}=t_I^n
x_\alpha(a) t_I^{-n}\in B^{(n)}$. By (\ref{txt}) we have $t_I^n u
t_I^{-n}=x_\alpha (\alpha(t_I^n)a)$, and by (\ref{TA}) we have
$$
w(\alpha(t_I^n))=-<\nu(t_I^n),\al>=n<\nu(t_I),-\alpha>=n<\sum_{j\in
\Delta\backslash I}\nu(t_j),-\alpha> .
$$
Now since $\alpha \in \Psi_I^-$, we have $-\alpha=\sum_{i\in
\Delta}n_i\alpha_i$ with $n_{i}\geq 0$ for any $i\in \Delta$ and
$n_j\geq 1$ for at least one $j\in \Delta\backslash I$. Hence
$$
w(\alpha(t_I^n))=n<\sum_{j\in \Delta\backslash I}\omega_j,\sum_{i\in
\Delta}n_i\alpha_i>=n \sum_{j\in \Delta\backslash I}n_j \geq n.
$$
Therefore, we have $t_I^n u t_I^{-n}=x_\alpha(\alpha(t_I^n)a)$ with
$\omega(\alpha(t_I^n)a)=(\alpha(t_I^n))+\omega(a)\geq n+1$. We have
just proved $t_I^n U_{\alpha,1} t_I^{-n} \subseteq
U_{\alpha,n+1}\subseteq B^{(n)}$. \qed

\begin{cor}
 Any compact open subset in $G/P_I$ can be written, for any $n\geq 0$ big enough,
 as a finite disjoint union of subsets of the form $bt_I^nB_IP_I/P_I$ with $b\in G(O)$.
\end{cor}
\proof See also \cite[prop. 8 and 8'.ii]{Schneider1}. Notice that we
have
$$
t_I^nB_IP_I=t_I^n(B\cap
U_I^-)t_I^{-n}P_I=\coprod_{x}xt_I^{n+1}(B\cap U_I^-)t^{-n-1}P_I
=\coprod_x xt_I^{n+1}B_IP_I
$$
where $x$ runs through the left cosets of $t_I^{n+1}(B\cap
U_I^-)t_I^{-n-1}$ in $t_I^{n}(B\cap U_I^-)t_I^{-n}$. \qed

As a consequence of this, we get the following result which implies
that the generalized Steinberg representations of $G$ are cyclic.

\begin{theo}\label{cyclicI}
The $M[G]$-module $C^{^\infty}(G/P_{I},\,M)$ is generated by the
element $\chi_{B_{I}P_{I}}$.
\end{theo}

\subsection{Generalized Steinberg representations and parahoric subgroups}

Let $G$ be the group of $K$-rational points of a semi-simple linear
algebraic group defined and split over $K$, of rank $l$ and of
adjoint type. Let $M$ be an abelian group.

Let $X$ be a totally discontinuous topological Space. We denote by
 $C_c^{^\infty}(X,\,M)$ the space of $M$-valued compactly supported locally
 constant functions on X.

Notice that, for any $I\subseteq \Delta$, there is an action of $G$
 on the space  $C_{c}^{^\infty}(G/B_I,\, M)$ induced by its left action on $G/B_I$.

In order to interpret the generalized Steinberg representations in
terms of parahoric subgroups, we start by considering the homomorphsim of $M[G]$-modules :
$$
H : C_{c}^{^\infty}(G/B,\, M) \longrightarrow C^{^\infty}(G/P,\, M)
$$
defined by $H(\varphi)=\varphi \star
\chi_{BP}:=\displaystyle\sum_{g\in G/B}\varphi(g).g(\chi_{BP})$.
Since we have
$$
B_IP_I=B_I P= \left(\coprod_{g\in B_I/B}gB\right)P=\coprod_{g\in
B_I/ B} gBP,
$$
where the last union is disjoint by Proposition \ref{disjoint}, we
get
$$
H(\chi_{B_I})=\sum_{g\in G/B}\chi_{B_I}(g).g(\chi_{BP})=\sum_{g\in
B_I/B}g(\chi_{BP})=\chi_{B_IP}=\chi_{B_IP_I}.
$$
Thus $H$ induces an $M[G]$-homomorphism which by Theorem
\ref{cyclicI} is surjective :
$$
H_I : C_{c}^{^\infty}(G/B_I,\,M) \longrightarrow
C^{^\infty}(G/P_I,\,M).
$$
What essentially remains to do is to compute the kernel of this
homomorphism. For this point we need some technical lemma. First,
recall that the convex closure $\textrm{cl}(\Omega)$ of a subset
$\Omega$ of $A$ is the intersection of all the closed half spaces in
$A$ that contain $\Omega$, i.e. we have
$$
\textrm{cl}(\Omega)= \bigcap_{{\al\in \Phi,r\in {\mathbb
Z}}\atop{\overline{H^+}_{\al,r}\supseteq \Omega}}
\overline{H^+}_{\al,r}
$$
where $\overline{H^+}_{\al,r}=\{x\in A;\, \langle x,\alpha \rangle
+r \geq 0 \}$.
\begin{lem} Let $I\subseteq \Delta$. We have
$$
\textrm{cl}(C\cup {\mathfrak C}_{I})=\bigcap_{\al\in \Phi_0}\overline{H^+}_{\al,0}\cap \bigcap_{\al\in
\Phi_I^-}\overline{H^+}_{\al,1}.
$$
\end{lem}
\proof It is obvious that $C \cup {\mathfrak C}_I$ is contained in any half plane involved
in the right hand side of the equality. So is the convex closure $\textrm{cl}(C\cup {\mathfrak C}_{I})$.
To prove the other inclusion, take any half plane
$\overline{H^+}_{\al,r}$, $\al\in \Phi$ and $r\in {\mathbb Z}$,
that contains $C\cup {\mathfrak C}_I$, it is clear that it should
contain the origin $0$ which is a vertex of $C$, this means that we must
have $r\geq 0$. Now, in case $\al \in \Phi^+$, we have
$\overline{H^+}_{\al,r}\supseteq \overline{H^+}_{\al,0}$ and the intersection
of all the half planes $\overline{H^+}_{\al,0}$ is the same as the intersection of
those with $\al\in \Phi_0$ and this is the topological closure of the sector
$\mathfrak C$. In case $\al\in \Phi^-$, if $\al \in \Phi^-_I$, so
$\al=\sum_{i\in I}m_i\al_i$ with some $m_j< 0$, the vertex $\om_j/n_j$
of $C$ should belong to $\overline{H^+}_{\al,r}$, so we should have
$\langle \om_j/n_j,\al \rangle +r \geq 0$, thus $r\geq -m_j/n_j >0$,
therefore we must have $r\geq 1$, thus $\overline{H^+}_{\al,r}\supseteq \overline{H^+}_{\al,1}$.
In case $\al\in \Phi^-\backslash\Phi_I$, we then have
$\al=\sum_{i\in \Delta}m_i\al_i$ with $m_j<0$ for some $j\in \Delta\backslash I$.
For any positive integer $n$, since $n\omega_j$ is a vertex of ${\mathfrak C}_I$ it
belongs to $\overline{H^+}_{\alpha, r}$, therefore $\langle n\omega_j,\alpha\rangle +r\geq 0$;
thus $r\geq -nm_j$ for any positive integer $n$, and this is impossible. \qed

\begin{theo}\label{kernelI}
Let $I\subseteq \Delta$ be such that any positive root
$\al=\sum_{i\in I}m_i\al_i$ in the root subsystem $\Phi_I$ satisfy
$m_i\leq 1$ for any $i\in I$. The kernel of the surjective map
$$
H_I : C_{c}^{^\infty}(G/B_I,\,M) \longrightarrow
C^{^\infty}(G/P_I,\,M)
$$
is the $M[G]$-submodule $R_I$ of $C_c^{^\infty}(G/B_I,M)$ generated
by the functions $\chi_{Bt_jB_I}-\chi_{B_I}$, $1\leq j \leq l$.
\end{theo}
\proof There is no need to give all the proof here since it is the same as given for the group $GL_{l+1}$ in \cite[Prop. 15]{Schneider1}. However, the following argument needs to be checked out in the more general case
$$
t_{j}^{-1}(B\cap P_{I})t_{j}\subseteq B_{I}\cap P_{I} \quad \textrm{ for any } 1\leq j\leq l.
$$
This argument which is in [loc.cit., Lemma 16] as an exercise is valid if and only if any positive root $\al=\sum_{i\in I}m_i\al_i$ in the root subsystem $\Phi_I$ satisfy
$m_i\leq 1$ for any $i\in I$. Indeed, when $j\in \Delta\backslash I$, we have $t_j\in T_I^{++}$, so by (\ref{TIBIPI}) we get
$$
t_j^{-1}(B\cap P_I)t_j \subseteq t_j^{-1}(B_I\cap P_I)t_j \subseteq
B_I \cap P_I.
$$
Now, using the fact that the stabilizer of a set and of its convex closure is the same, we have
$$
B\cap P_{I}=B\cap B_{I} \cap P_{I}={\cal P}_{C} \cap {\cal
P}_{{\mathfrak C}_{I}}={\mathcal P}_{C\cup {\mathfrak C}_I}={\mathcal P}_{\textrm{cl}(C\cup {\mathfrak C}_I)}.
$$
Therefore, our argument is equivalent to the following
$$
{\mathcal P}_{\textrm{cl}(C\cup {\mathfrak C}_I)}\subseteq t_j{\mathcal
P}_{{\mathfrak C}_I}t_j^{-1} \qquad \textrm{for any }j\in I
$$
or to
$$
\omega_j + {\mathfrak C}_I \subseteq \textrm{cl}(C \cup {\mathfrak C}_I) \qquad \textrm{for any }j\in I.
$$
It remains to notice that for $\alpha=\sum_{i\in I}m_i\alpha_i \in\Phi_I^-$ and $x\in {\mathfrak C}_I$, we have
$$
\langle \omega_j + x,\alpha\rangle +1= m_j +1.
$$
So $\omega_j +x \in \overline{H^+}_{\alpha,1}$ if and only if $m_j\geq -1$. \qed

As consequences, we get the generalized Steinberg representations of $G$ in terms of the parahoric subgroups :

\begin{cor}\label{stI} Let $I\subseteq \Delta$ be as in Theorem \ref{kernelI} above. We have
an isomorphism of $M[G]$-modules :
$$
\textrm{St}_I(M)\cong \frac{C_c^{^\infty}(G/B_I,M)}{R_I+\sum_{i\in
\Delta\backslash I} C_c^{^\infty}(G/B_{I\cup\{i\}},M)}
$$
where $R_I$ is the $M[G]$-submodule of $C_c^{^\infty}(G/B_I,M)$
generated by the functions $\chi_{Bt_jB_I}-\chi_{B_I}$, $1\leq j
\leq l$.
\end{cor}
\proof Apply the theorem above. \qed

\begin{cor}\label{st} We have an isomorphism of $M[G]$-modules :
$$
St(M)\cong \frac{C_c^{^\infty}(G/B,M)}{R+\sum_{i\in \Delta}
C_c^{^\infty}(G/B_{\{i\}},M)}
$$
where $R$ is the $M[G]$-submodule of $C_c^{^\infty}(G/B,M)$
generated by the functions $\chi_{Bt_jB}-\chi_{B}$, $1\leq j \leq
l$.
\end{cor}
\proof Just take $I=\emptyset$ in the corollary above.

\bigskip
Y. A\"IT AMRANE, \\
Facult\'e de Math\'ematiques,\\
USTHB BP32, Al Alia,\\
16111 Bab-Ezzouar, ALGER.\\
e-mail adress : yaitamrane@usthb.dz or amrane@math.cnrs.fr

\end{document}